\documentclass[12pt]{gtart}
\usepackage{amsmath}
\usepackage{amssymb}
\usepackage{amsthm}
\usepackage{pinlabel}

\newtheorem{theorem}{Theorem}

\theoremstyle{definition}

\newtheorem{remark}[theorem]{Remark}

\def\R{\mathbb R}

\newcommand{\Diff}{\mathop{\rm Diff}\nolimits}

\begin{document}

\title{Reconstructing $4$-manifolds from Morse $2$-functions}
\authors{ David T. Gay, Robion Kirby
\footnote{This work was partially supported by a grant from the Simons Foundation (\#210381 to David Gay). The second author was partially supported by NSF grant DMS-0838703.}}
\address{Euclid Lab, 500 Willow St, Athens, GA 30601\\
Department of Mathematics,University of Georgia, Athens, GA 30602}
\secondaddress{University of California, Berkeley, CA 94720} 
\email{d.gay@euclidlab.org}
\secondemail{kirby@math.berkeley.edu}
\begin{abstract}
Given a Morse $2$-function $f \co X^4 \to S^2$, we give minimal conditions on the fold curves and fibers so that $X^4$ and $f$ can be reconstructed from a certain combinatorial diagram attached to $S^2$.  Additional remarks are made in other dimensions.
\end{abstract}
\primaryclass{57M50}
\secondaryclass{57R45, 57R65}
\keywords{broken fibration, Morse function, Morse 2-function, Cerf theory, 4-manifold}

\maketitle

Let $X^n$ and $\Sigma^2$ be smooth, oriented, closed, connected
manifolds. A {\it Morse $2$-function} $f \co X^n \to \Sigma^2$ is a
smooth function such that each point $x_0 \in X^n$ has a coordinate
chart $\R \times \R^{n-1}$ and $f(x_0)$ has a coordinate chart $\R \times
\R$ for which $f(x) = f(t,y) = (t,f_t(y))$ where $f_t \co \R^{n-1} \to \R$ is
a generic $1$-parameter family of smooth functions with non-degenerate
critical points and births and deaths~\cite{gk1}. The image of the set of critical points of $f_t$ is
often called a {\it Cerf graphic}~\cite{cerf}.  The set of points in
$X$ for which the rank of $Df$ is only one forms a smoothly embedded
$1$-manifold $Z$ whose image in $\Sigma$ is an immersed $1$-manifold
$Z'$ (with cusps) of critical values of $f$.  $Z'$ is called the
set of {\it fold curves} in $\Sigma$.  Under appropriate conditions on the homotopy class of the Morse $2$-function $f$, in particular when $\Sigma=S^2$, $f$ is homotopic
to a Morse $2$-function for which all fibers are connected and all critical points are indefinite~\cite{gk1}. That is, in local
coordinates, $f(t,y) =  (t, -y_1^2 \pm y_2^2 \pm \dots \pm y_{n-2}^2 +
y_{n-1}^2)$. From now on we assume all Morse $2$-functions are indefinite and fiber-connected.

\begin{figure}
\labellist
\small\hair 2pt
\pinlabel $a$ [b] at 241 30
\pinlabel $b$ [t] at 241 94
\pinlabel $a$ [l] at 254 62
\pinlabel $b$ [r] at 228 62
\pinlabel $a$ [bl] at 215 20
\pinlabel $a$ [tl] at 215 14
\pinlabel $b$ [br] at 213 106
\pinlabel $b$ [bl] at 216 106

\endlabellist
\centering

 \includegraphics[width=5in]{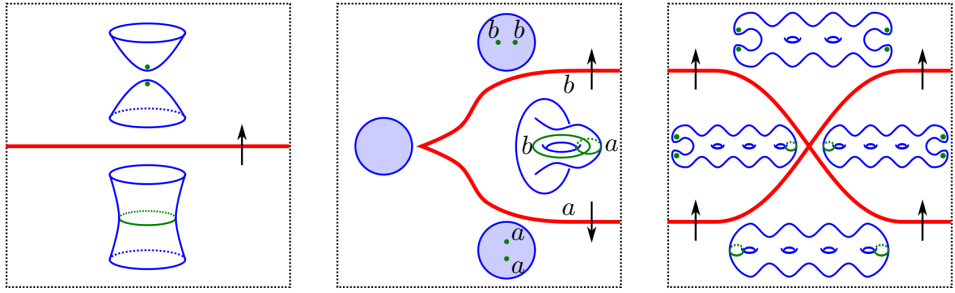}
 \caption{\label{F:LocalModels} Local models for indefinite fiber-connected Morse $2$--functions on $4$--manifolds. Folds are in red, fibers (or local parts of fibers) are in blue, attaching $S^0$'s and $S^1$'s are in green. Arrows transverse to folds indicate the direction of decreasing genus, i.e. the direction in which a fold has index $2$.}
\end{figure}

In the case $n=4$ which we focus on in this paper, all the fibers
(preimages of points under $f$) are oriented surfaces, and when
crossing a fold curve in $\Sigma$ the genus changes by one. On the
higher genus side is a circle in the fiber to which a $2$-handle is
attached (we call this the {\em attaching circle} for the fold) and on the lower genus side are a pair of points to which a
$1$-handle is attached;  these are the ascending and descending
spheres of the critical point of $f_t$ at a fold curve.  At a cusp,
the pair of circles, one from each side of the cusp, must meet in
exactly one point in the fiber;  at a crossing of fold curves, they
must be disjoint (see Figure~\ref{F:LocalModels}).

\begin{figure}
\labellist
\small\hair 2pt
 \pinlabel $A_1$ [t] at 22 24
 \pinlabel $B_1$ [l] at 31 36
 \pinlabel $A_2$ [t] at 59 21
 \pinlabel $B_2$ [l] at 66 37
 \pinlabel $A_3$ [t] at 96 22
 \pinlabel $B_3$ [l] at 106 36
 \pinlabel $A_g$ [t] at 181 25
 \pinlabel $B_g$ [r] at 166 34

\endlabellist
\centering

 \includegraphics[width=5in]{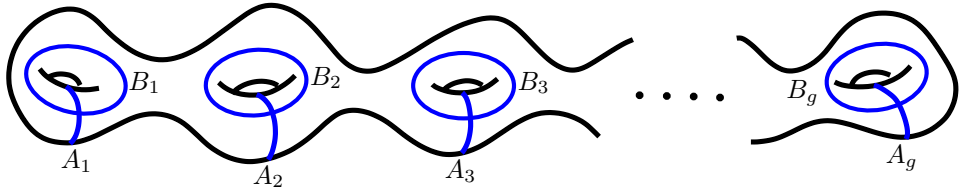}
 \caption{\label{F:StandardBasis} The standard genus $g$ surface and the standard basis $A_1, B_1, \ldots, A_g, B_g$.}
\end{figure}

\begin{theorem} \label{T:BasicData}
Suppose that all the fibers of a given Morse $2$--function $f \co X \to \Sigma$ have genus $>1$, and
all the regions in $\Sigma$ bounded by fold curves are simply connected. Then the following data suffice to reconstruct $X$ and $f$ up to diffeomorphism:
\begin{enumerate}
 \item The {\em fold graph} $\Gamma$, the image of the fold curves in $\Sigma$. This is a co-oriented graph with vertices of valence $2$ and $4$ dividing $\Sigma$ into polygonal regions; the co-orientation is in the direction of decreasing genus. When $\partial \Sigma \neq \emptyset$, the vertices should be in the interior and edges should be transverse to the boundary.
 \item The {\em standard fiber} in each region $R$, a drawing of a standard genus $g_R$ surface $F_R$; this is the fiber over a point in the interior of $R$.
 \item The {\em attaching circle} for each outward oriented edge $e$ of a region $R$, a simple closed curve $C_e$ drawn on $F_R$; this is the attaching circle for the $2$--handle attached when crossing $e$.
 \item The {\em gluing data} for each outward oriented edge $e$ of a region $R$, a basis of $2(g_R-1)$ simple closed curves $A_{(1,e)}, B_{(1,e)}, \ldots, A_{((g_R-1),e)}, B_{((g_R-1),e)}$ disjoint from $C_e$, with $A_i \cap B_i$ a single transverse point of intersection and $(A_i \cup B_i) \cap (A_j \cup B_j) = \emptyset$ when $i \neq j$; this specifies that the genus $g_R-1$ surface obtained from surgery along $C_e$ should be identified with the standard fiber in the adjacent region across edge $e$ so that the basis $A_{(1,e)}, B_{(1,e)}, \ldots, A_{((g_R-1),e)}, B_{((g_R-1),e)}$ maps to the standard basis for the standard genus $g_R-1$ fiber. (Figure~\ref{F:StandardBasis} illustrates this standard basis.)
\end{enumerate}
\end{theorem}

\begin{remark}
Without item~4 in the list above (the gluing data) and with $\Sigma=S^2$, this is Theorem~3 of \cite{gk2}. {\em That theorem is in fact false}. We will present explicit counterexamples showing that, in general, one does need to record {\em at least some} of the gluing data. However, in the next theorem we will see that one should expect many edges across which one does {\em not} need to record the gluings.
\end{remark}

\begin{remark}
 If one draws an arbitrary co-oriented $2$-- and $4$--valent graph in a surface $\Sigma$ with data as in the theorem, there are many obvious constraints that would prevent such data from describing a closed $4$--manifold with a Morse $2$--function $f$. Obvious constraints include that the genus should drop by one across each edge, that co-orientations should be consistent at $4$--valent vertices, and that attaching circles at crossings should be disjoint and at cusps should intersect once. There are also much less obvious constraints on both attaching circles and gluing maps related to the fact that, after reconstructing $X$ over a neighborhood of a dual graph, we need to be able to fill in the inverse images of disks containing the cusps and crossings.
\end{remark}

\begin{remark}
 It is almost immediate that the data in Theorem~\ref{T:BasicData} will reconstruct $(X,f)$ over the complement of disk neighborhoods of the cusps and crossings, and thus the real content of the theorem is that, once this reconstruction is done, the extension over the cusps and crossings is unique if it exists.
\end{remark}

\begin{figure}
\labellist
\small\hair 2pt
\endlabellist
\centering

 \includegraphics{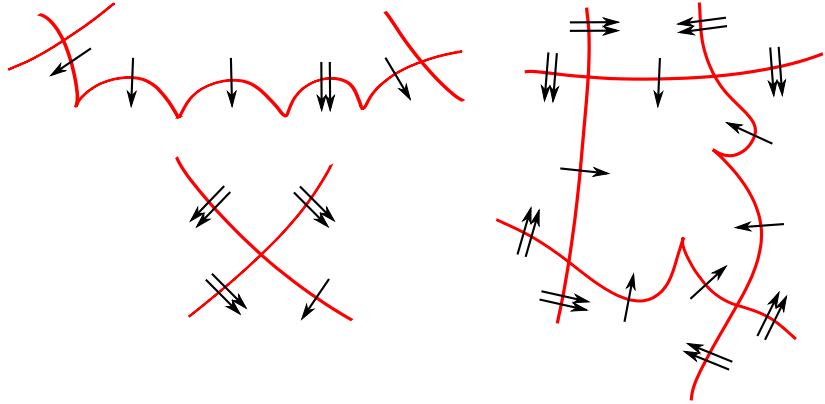}
 \caption{\label{F:NoGluingData} Cases in which we do not need to record all the gluing data. The edges across which the gluing data needs to be recorded are labelled with double transverse arrows; the single transverse arrows indicate edges for which the gluing data can be recovered from the gluing data along the double arrows.}
\end{figure}

\begin{theorem} \label{T:NoGluingData}
The following are situations under which we only need to record the gluing data across a subset of the edges of the fold graph $\Gamma$. See Figure~\ref{F:NoGluingData}.
\begin{enumerate}
 \item Define a {\em chain of edges} to be a sequence of edges connected by cusps, and thus all adjacent to the same two regions. For any chain of edges, the gluing data across one edge in the chain determines the gluing data across all other edges in the chain.
 \item At a crossing involving one region of genus $g+2$, two of genus $g+1$ and one of genus $g$, consider the four edges meeting at that crossing. The gluing data across the two edges between genus $g+2$ and genus $g+1$ together with the gluing data across one of the edges between genus $g+1$ and genus $g$ determines the gluing data across the remaining edge (between genus $g+1$ and genus $g$).
 \item If $R$ is a locally genus minimizing region, so that all the edges of $R$ are co-oriented inwards, then we do not need to record any gluing data across any edges of $R$. In other words, the gluing data across all the other edges meeting $R$ at vertices (crossings) determines the gluing data across the edges of $R$ up to an automorphism of the fiber $F_R$, and changing the gluing data across the edges of $R$ by such an automorphism does not change $(X,f)$.
\end{enumerate}
\end{theorem}

\begin{remark}
 Theorem~\ref{T:NoGluingData} implies that in many cases, one only needs to record the gluing data along a rooted maximal tree $T$ in the oriented dual graph $\Gamma^*$ to $\Gamma$. The condition we need is that, in the process of adding edges to turn $T$ into $\Gamma^*$, whenever one extends across a quadrilateral (which corresponds to a crossing in $\Gamma$), the edge being added is one of the two lower-genus edges of the quadrilateral. For many graphs such a maximal tree will not exist.
\end{remark}

To warm up, consider a traditional Morse function $g \co X^n \to \R$.  The
analogues of fold curves are critical values, and the analogues of
fibers are level sets, and the analogues of circles and pairs of
points are the ascending and descending spheres.  This is not enough
information to reconstruct $X^n$ in general.  One needs to know either 
\begin{enumerate}
\item how to identify the level set just above a critical value with
  the level set just below the next higher critical value (this can be
  done with a choice of metric which gives a gradient flow), or
\item how to see all the attaching maps at once in the boundary (=
  level set) of the $0$-handle, for all the handles (including even
  the last $n$-handle, as can be seen by considering exotic smooth
  structures on spheres).
\end{enumerate}

In the special case of dimension four, $X^4$ can be
determined by drawing the attaching maps of  $1$-handles as pairs of
points in $S^3$, the $2$-handles as framed links in $S^3$ which may go
over $1$-handles, and then using the homotopy type of $X^4$ to
determine how to add the $3$-handles and $4$-handle (see \cite{lp}).

When considering a circle valued Morse function, $g \co X^n \to S^1$, even
the information of (1) and (2) is not enough for $g$ may be a bundle
map with no critical values and hence nothing to determine the
monodromy of the bundle.

Right off we see that a non-trivial circle bundle $F^2 \to M^3 \to
S^1$ produces a Morse $2$-function $M^3 \times S^1 \to S^1 \times S^1$
with no fold curves and therefore no information to distinguish $M^3
\times S^1 $ from, say, $F \times S^1 \times S^1$.

Traditionally we define bundles $p \co X^n \to \Sigma^2$ by giving an open
cover $\{U_{\alpha}\}$ of $\Sigma$ and local trivializations
$h_{\alpha} \co p^{-1}(U_{\alpha}) \to U_{\alpha} \times F^{n-2}$ together
with clutching functions $\phi_{\alpha\beta} \co h_{\alpha}p^{-1}(x)
\to h_{\beta}(p^{-1}(x))$ for $x \in U_{\alpha} \cap U_{\beta}$, which
satisfy the cocycle condition, $\phi_{\alpha\beta}\phi_{\beta\gamma} =
\phi_{\alpha\gamma}$.

This definition can be extended to what might be called a broken
bundle $p \co S^n \to \Sigma^2$ in which $\Sigma$ has the usual fold
curves with their data of ascending and descending spheres specified
by giving an open cover $\{U_{\alpha}\}$ of $\Sigma$, and local
trivializations of the follow four kinds:

\begin{enumerate}

\item $h \co p^{-1}(U_{\alpha}) \to U_{\alpha} \times F^{n-2}$,

\item $h \co p^{-1}(U_{\alpha}) \to B$ where $B$ is an interval cross an
  $(n-1)$-dimensional bordism with one critical point,

\item $p^{-1}(U_{\alpha})$ is the obvious preimage of a birth or
  death of a canceling pair of critical points,

\item $p^{-1}(U_{\alpha})$ is the obvious preimage of a crossing of
  fold curves.

\end{enumerate}

In addition, clutching functions must be given for the traditional case
with no folds involved, or for gluing along a fold curve.  Again the
cocycle condition must be satisfied, and then the broken bundle, hence
$X^n$, is well defined.

With these remarks behind us, we now give the proof of Theorem~\ref{T:BasicData}.

\begin{proof}
Clearly the fold graph $\Gamma$ and the standard fibers determine $X$ and $f$ over a disk in the interior of each region, the blue disks in Figure~\ref{F:BasicExample}. Enlarge each of these disks by sending out fingers (green in Figure~\ref{F:BasicExample}) which intersect each outward-pointing edge along an arc in the interior of the edge, but are still disjoint from the neighboring disks; the attaching circles then determine $X$ and $f$ over these enlarged disks because, in $X$, adding each finger corresponds to attaching $I$ cross a $3$--dimensional $2$--handle cobordism for each edge, with the $2$--handle attached along the given attaching circle. 
\begin{figure}
\labellist
\small\hair 2pt

\endlabellist
\centering

 \includegraphics{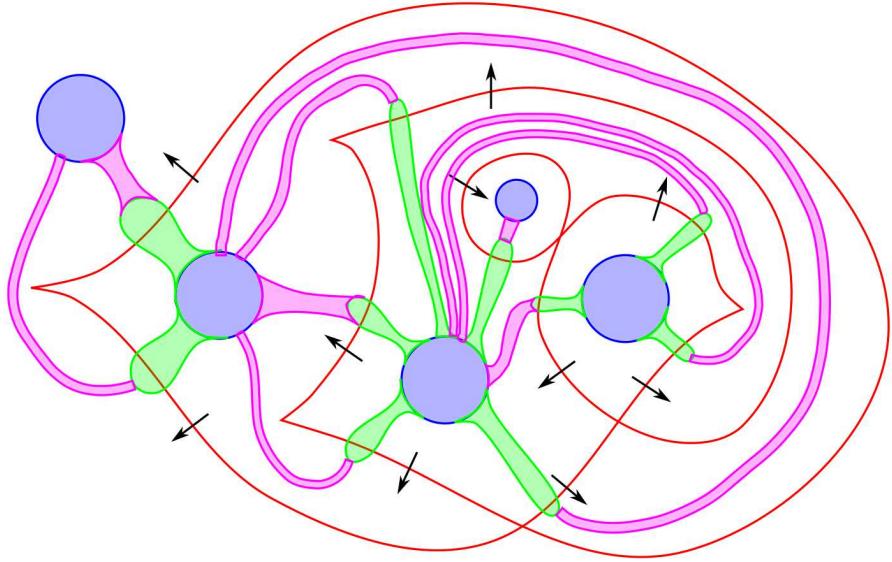}
 \caption{\label{F:BasicExample} Using the basic data to build $X$, first over blue disks, then green fingers, then pink bridges.}
\end{figure} 
The gluing data are then used to attach the tips of the fingers to the adjacent disks, along the pink bridges in Figure~\ref{F:BasicExample}: Consider two regions $R$ and $R'$ meeting along an edge $e$ co-oriented from $R$ to $R'$. Thus $g_{R'} = g_R - 1$. Let $F_{R,e}$ be the result of surgering $F_R$ along $C_e$. The attaching data $A_1, B_1, \ldots, A_{g_R-1}, B_{g_R-1}$ form a basis for $F_{R,e}$. Meanwhile, $F_{R'}$ is drawn in a standard way which implies a standard basis, and the gluing map is then constructed to send the given basis of $F_{R,e}$ to the standard basis on $F_{R'}$.

Thus the full data determines $(X,f)$ over the complement of disk neighborhoods of the cusps and crossings. We now need to show that the extensions over the cusps and crossings are unique. We know that the extensions exist by the assumption that this data comes from a given $(X,f)$. 

We will first show the uniqueness of three kinds of extensions across disks, where the fibration in each case is given over half of the boundary, as in Figure~\ref{F:3Extensions}. 
\begin{figure}
\labellist
\small\hair 2pt
 \pinlabel $(a)$ [t] at 62 0
 \pinlabel $(b)$ [t] at 202 0
 \pinlabel $(c)$ [t] at 342 0
 \pinlabel $p$ [r] at 309 60
 \pinlabel $q$ [l] at 379 60
 \pinlabel $P$ [r] at 312 92
 \pinlabel $Q$ [l] at 375 92
 \pinlabel $A$ [t] at 344 32

\endlabellist
\centering

 \includegraphics[width=5in]{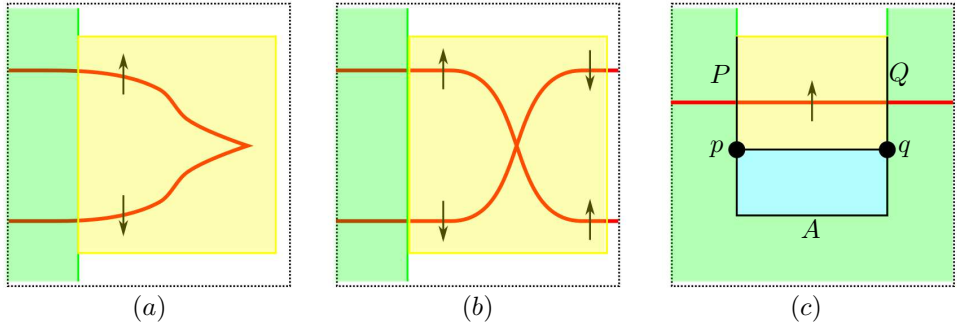}
 \caption{\label{F:3Extensions} Three cases in which an extension across a disk is unique. In the first two cases we are extending from the green region to the yellow region. In the last case we are extending first from the green to the yellow and then to the blue.}
\end{figure} 
In the three cases in Figure~\ref{F:3Extensions} we have chosen local coordinates so that we have Cerf diagrams. 

In Figure~\ref{F:3Extensions}(a) and~(b) we are extending an ordinary Morse function $g \co M^3 \to [0,1]$ to a stable $1$--parameter family $g_t \co M \to [0,1]$ with $g_0 = g$, with a prescribed Cerf graphic for $g_t$, either a cusp or a crossing. In the cusp case, the uniqueness assertion is then just the assertion that, if two critical points of index $1$ and $2$, with adjacent critical values, can be cancelled then there is a unique stable $1$-parameter family of Morse functions achieving this cancellation. In the crossing case, the direction in which we are extending is such that we start with a Morse function with an index $1$ critical point below an index $2$ critical point, and we know that the descending manifold for the index $2$ critical point is disjoint from the ascending manifold for the index $1$ point, and then we extend to a $1$--parameter family in which the index $1$ point rises above the index $2$ point. Up to diffeomorphism this is unique because we do not have any handle slides between the two critical points to worry about. (By comparison, if we were extending in a direction in which both critical points had the same index, there would be many different extensions depending on whether we put handle slides in before the critical point crossing.)

In Figure~\ref{F:3Extensions}(c), we have illustrated the extension in two steps, first from green to yellow and then to blue. In the first step (yellow), we need to choose a diffeomorphism $\phi$ from the fiber $F_p$ at $p$ to the fiber $F_q$ at $q$, sending the attaching circle $C_p \subset F_p$ for the $2$--handle cobordism $H_p$ sitting over the arc $P$ to the attaching circle $C_q \subset F_q$ for the $2$--handle cobordism $H_q$ sitting over the arc $Q$. The results of the first step only depend upon $\phi$ up to isotopy through diffeomorphisms sending $C_p$ to $C_q$. After the first step, the monodromy around the enclosed disk (blue) should be trivial. Thus $\phi$ is determined up to isotopy (not necessarily rel. $C_p$ and $C_q$) by the existing identification between $F_p$ and $F_q$ given by going around the lower side of the blue disk, the arc $A$. (This uses the fact that $\pi_1(\mathop{\Diff}_0\Sigma_g) =0$ when $g>1$; see~\cite{earleeells}.) Now, since all fibers have genus larger than $1$, maps sending $C_p$ to $C_q$ which are isotopic are actually isotopic rel. $C_p$ and $C_q$, and thus the isotopy class of $\phi$ relative to the attaching curves is determined.

Given the uniqueness of the three extensions above, the uniqueness of extension across a cusp or crossing from the full boundary of a disk to the interior of the disk follows as in Figure~\ref{F:FullExtensions}, where the last extension across a disk containing no folds is again unique because $\pi_1(\mathop{\Diff}_0\Sigma_g) =0$. It is important that we apply all the extensions from Figure~\ref{F:3Extensions} in the correct direction, i.e. in the direction of decreasing genus.
\begin{figure}
\labellist
\small\hair 2pt
\endlabellist
\centering

 \includegraphics{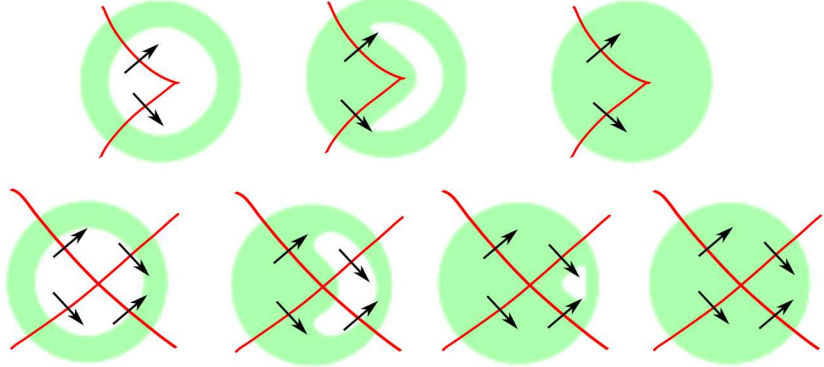}
 \caption{\label{F:FullExtensions}Full extensions across cusps and crossings.}
\end{figure}

\end{proof}

\begin{proof}[Proof of Theorem~\ref{T:NoGluingData}]
 The first two cases (chains of edges and crossings) both follow immediately from the uniqueness of the three extensions across disks shown in Figure~\ref{F:3Extensions}. The case of a locally genus minimizing region is also immediate: It is clear from the first two cases that the gluing data across one edge of $R$ determines the gluing data across all other edges. Since there are no outward co-oriented edges of $R$, changing this gluing data by an automorphism of $F_R$ changes no other data in the diagram, and any two gluing data across one fixed edge of $R$ are related by such an automorphism.
\end{proof}

We next present two interesting counterexamples to our original false claim (Theorem~3 of~\cite{gk2}), that in the case of maps to $S^2$ the gluing data is unnecessary. The first is a local picture that can be embedded in a $S^2$--valued Morse $2$--function on some closed $4$--manifold, and the second is a more global construction.

\begin{figure}
\labellist
\small\hair 2pt
\pinlabel $a$ [t] at 33 61
\pinlabel $b$ [t] at 83 68
\pinlabel $c$ [tl] at 213 68
\pinlabel $d$ [t] at 299 56
\pinlabel $u$ [b] at 40 117
\pinlabel $x$ [br] at 127 109
\pinlabel $y$ [b] at 157 112
\pinlabel $z$ [bl] at 201 109
\pinlabel $v$ [b] at 272 111
\pinlabel $a$ [b] at 121 39
\pinlabel $b=d$ [t] at 127 3
\pinlabel $c$ [t] at 183 3
\pinlabel $u$ [l] at 35 83
\pinlabel $x=z$ [t] at 156 70
\pinlabel $y$ [l] at 171 85
\pinlabel $v$ [l] at 275 78
\pinlabel $A$ [tr] at 167 5
\pinlabel $B$ [l] at 201 27
\endlabellist
\centering

 \includegraphics{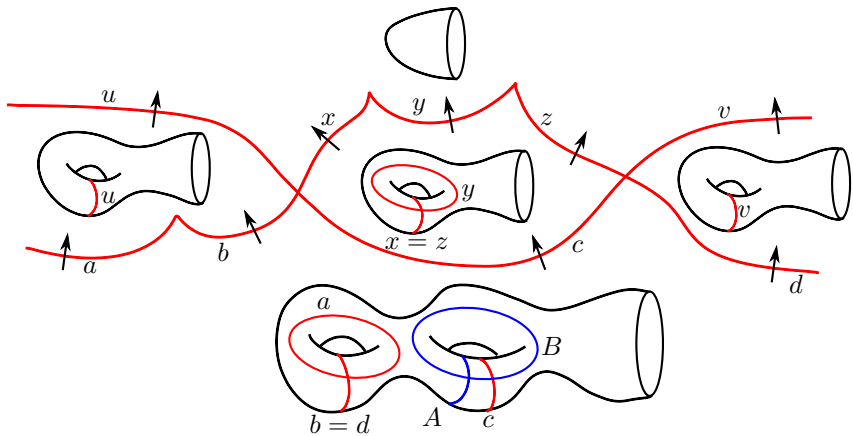}
 \caption{\label{F:CounterExLocal} The local counterexample. Curves in red are attaching circles. The curves in blue labelled $A$ and $B$ are the gluing data for both edges $a$ and $d$. In other words, $A_a = A_d = A$ and $B_a = B_d = B$. No gluing data is shown for other edges.}
\end{figure}

In Figure~\ref{F:CounterExLocal} we show the fold graph $\Gamma$ in a disk, with the attaching circles indicated, and the gluing data specified around the boundary. To avoid excess subscripts, the attaching circles are labelled with the same label as the edges. The fiber is drawn as a surface with boundary, the understanding being that this can glued onto a surface of any genus and the gluing data extended by the identity map. Recall that the gluing data across an edge $e$ is a basis $A_{(1,e)}, B_{(1,e)}, \ldots, A_{((g_R-1),e)}, B_{((g_R-1),e)}$ disjoint from the attaching circle for that edge. Thus for the upper edges, going from genus $1$ to genus $0$, there is nothing to draw, and for the lower edges, we should draw one pair $(A_e, B_e)$ for each edge $e$. The gluing data then means that, after surgering the genus $2$ surface along the attaching circle for edge $e$, we glue the resulting genus $1$ surface to the standard genus $1$ surface sending $(A_e,B_e)$ to the standard basis (see Figure~\ref{F:StandardBasis}).

\begin{figure}
\labellist
\small\hair 2pt
\pinlabel $a$ [t] at 33 124
\pinlabel $b$ [t] at 83 131
\pinlabel $c$ [tl] at 213 131
\pinlabel $d$ [t] at 299 121
\pinlabel $u$ [b] at 40 180
\pinlabel $x$ [br] at 127 172
\pinlabel $y$ [b] at 157 175
\pinlabel $z$ [bl] at 201 172
\pinlabel $v$ [b] at 272 174
\pinlabel $y$ [l] at 171 148
\pinlabel $A_c$ [tr] at 12 72
\pinlabel $B_c$ [r] at 4 94
\pinlabel $A_c$ [tr] at 173 71
\pinlabel $B_c$ [bl] at 174 99
\pinlabel $y$ [l] at 55 26
\pinlabel $y$ [l] at 215 24
\pinlabel $c$ [t] at 83 65
\pinlabel $c$ [t] at 243 64
\pinlabel $A$ [tr] at 68 68
\pinlabel $B$ [l] at 102 89
\pinlabel $A$ [tr] at 228 67
\pinlabel $B$ [l] at 262 88
\endlabellist
\centering
\includegraphics{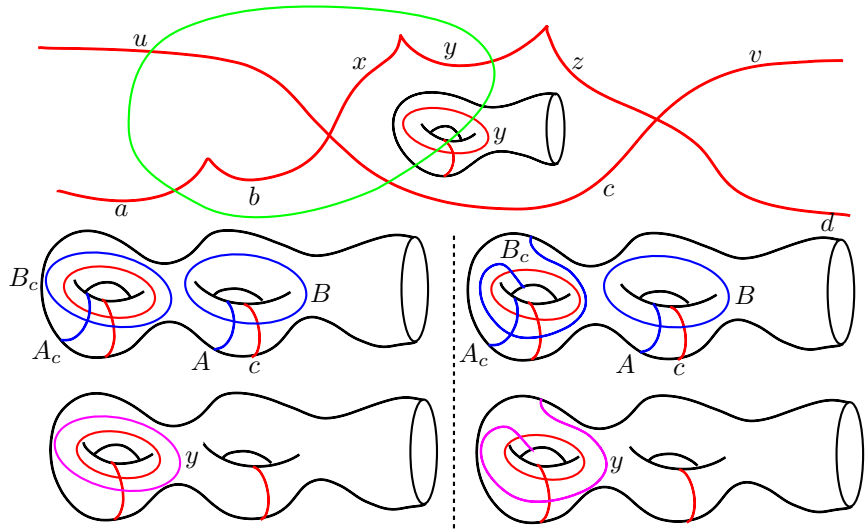}
 \caption{\label{F:CounterExLocalPt2} The local counterexample, part 2. Two choices for the gluing data $(A_c,B_c)$ across edge $c$ are indicated in blue, one to the left of the dotted line, the other to the right. The pullbacks by these two gluings of the attaching circle for edge $y$ are indicated in pink on the lower surfaces. The curves $A$ and $B$ are the gluing data across edges $a$, $b$ and $d$, i.e. $A=A_a=A_b=A_d$ and $B=B_a=B_b=B_d$. All unlabelled curves should be taken to be labelled as in Figure~\ref{F:CounterExLocal}.}
\end{figure}

The ambiguity lies in the gluing data across edge $c$ at the bottom, and in particular, the curve $B_c$ in the basis $(A_c,B_c)$. The pair $(A_c,B_c)$ should be a basis for what's left after one surgers along $c$, and thus should be disjoint from $c$. For example, $(A_c,B_c)$ can be either of the two examples illustrated in Figure~\ref{F:CounterExLocalPt2}. To see that the results are different, pull back the attaching circle for edge $y$ from the genus $1$ surface to the genus $2$ surface via the gluing map; in one case, $y$ becomes a circle that intersects $a$ once, transversely, while in the other case, $y$ is disjoint from $a$. Thus, if one considers the $S^1$--valued Morse function sitting over the green loop in Figure~\ref{F:CounterExLocalPt2}, either the critical points $a$ and $y$ cancel or they do not, and so the homology of the $3$--manifold sitting over the blue circle is different.

\begin{figure}
\labellist
\small\hair 2pt
\pinlabel $*$ [c] at 115 172
\pinlabel $a$ [t] at 115 137
\pinlabel $b$ [t] at 115 161
\pinlabel $c$ [b] at 115 186
\pinlabel $d$ [b] at 115 209
\pinlabel $b=d$ [bl] at 339 127
\pinlabel $a=c$ [bl] at 319 197
\pinlabel $\alpha$ [b] at 293 202
\pinlabel $\beta$ [br] at 302 126
\pinlabel $\gamma$ [tr] at 294 27
\pinlabel $\delta$ [br] at 292 55
\pinlabel $\rho$ [r] at 126 66 
\pinlabel $\sigma$ [r] at 112 39 
\pinlabel $\tau$ [r] at 90 17  
\endlabellist
\centering
 \includegraphics[width=4in]{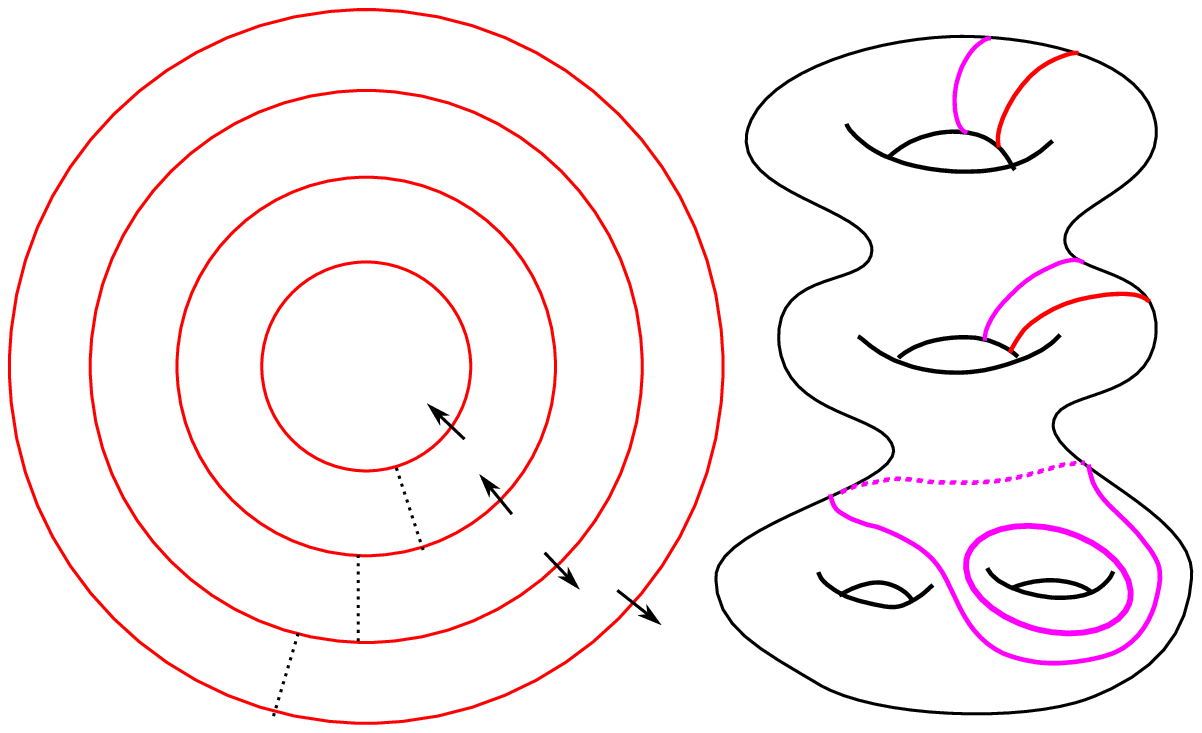}
 \caption{\label{F:CounterExGlobalPt1} The global counterexample. The genus $4$ surface drawn at right is the fiber at the point labelled $*$, with all attaching circles drawn in red. The pink curves are Dehn twist curves representing the monodromy across the dotted lines: $\rho = \alpha^k \gamma^m \delta^{-m}$, $\tau = \beta^l \gamma^m \delta^{-m}$, $\sigma = \alpha^k \beta^l \gamma^m \delta^{-m}$, for some $k,l,m \in \mathbb{Z}$. }
\end{figure}

\begin{figure}
\labellist
\small\hair 2pt
\endlabellist
\centering
 \includegraphics[width=2in]{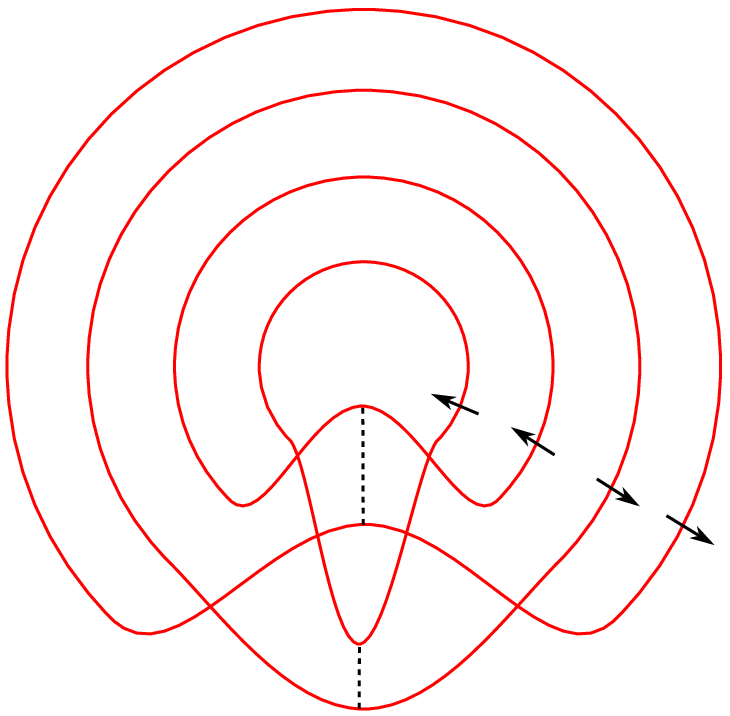}
 \caption{\label{F:CounterExGlobalPt2} The global counterexample, part 2.}
\end{figure}

In Figure~\ref{F:CounterExGlobalPt1} we show a Morse $2$--function over $S^2$ with annular regions between the folds, and we indicate the possibilities for monodromy around those annuli. We think of this monodromy as being a self-gluing of the fiber along the vertical dotted lines. Notice that we cannot choose any monodromy, because the monodromies should become trivial after surgering along the attaching circles for the edges, so that the disks can be filled in. We have indicated here a few of the many possible choices satisfying this condition. In Figure~\ref{F:CounterExGlobalPt2} we attempt to hide the monodromy by making all regions simply connected. The Reidemeister II type moves that we have done in the middle to cut the annuli into disks are legitimate because all the attaching circles involved are disjoint. Since the monodromy is now ``happening'' in the two small quadrilateral regions in the middle, it does not affect any of the attaching circles for the edges, and thus the attaching circles alone cannot possibly be sufficient to reconstruct the monodromy.

This last example leaves us with an interesting question: How can one formalize the notion of nontrivial monodromy when all regions and even the entire base are simply connected, without reference to a homotopic fibration with non simply connected regions?

\end{document}